\documentclass{article}

\usepackage{graphicx}

\usepackage{amsmath}
\usepackage{amssymb}
\usepackage{indentfirst}
\usepackage{enumerate}

\begin{document}

\title{Solving Poisson equations by the MN-curve approach }         
\author{Lin-Tian Luh\\Department of Data Science and Big Data Analytics, Providence University\\Shalu, Taichung, Taiwan\\
Email: ltluh@pu.edu.tw\\Fax: 886-4-26324653, Tel: 886-4-26328001 ext. 15126 }        
\date{\today}          
\maketitle
\hspace{-5.5mm}{\bf Abstract} In this paper we apply the newly born choice theory of the shape parameters contained in the smooth radial basis functions to solve Poisson equations. Some people complain that Luh's choice theory, based on harmonic analysis, is mathematically complicated and applies only to function interpolations. Here we aim at presenting an easily accessible approach to solving differential equations with the choice theory which proves to be successful, not only by its easy accessibility, but also by its striking accuracy and efficiency.

\section{Introduction}
Here we focus on the generalized multiquadrics
\begin{align}
\phi(x):=(-1)^{\lceil \beta/2\rceil}(c^{2}+\| x\|^{2})^{\beta/2},\quad \beta \in R\setminus 2N_{\geq 0},\quad c>
0,\quad x \in R^{n},
\end{align}
where $\lceil \beta/2\rceil$ denotes the smallest integer greater than or equal to $ \beta/2$. These are the most popular radial basis functions (RBFs) and frequently used in the collocation method of solving partial differential equations. The choice of the shape parameter $c$ contained in $\phi(x)$ has been obsessing experts in this field for decades and often leads to giving up this approach. Hitherto there is no theory about its optimal choice when dealing with PDEs. Although Luh's theory, called the c-theory by E. Kansa, can predict its optimal value almost exactly, it applies to function interpolations only and involves the complicated theory of harmonic analysis. Scientists, especially non-mathematicians, still do not know how to choose it when solving PDEs with RBFs.

As the inventor of the choice theory, the author knows that the theory applies to PDEs as well, maybe with a moderate search when necessary. The main reason is that collocation is in spirit a kind of interpolation. Moreover, Dirichlet conditions do offer interpolation points on the boundary. As can be seen in Luh \cite{Lu1,Lu2}, when $c$ is chosen according to the MN curves, the accuracy of the function approximation is incredibly good, both in theory and practice. It is not hard to imagine that the combination of the c-theory and collocation may lead to subversive results in the field of numerical PDEs.

In this paper we follow Kansa's route \cite{Ka1,Ka2,Ka3} to make collocation, but in a totally different way of choosing $c$. Basically, we discard the traditional trial-and-error search, and adopt the theoretically predicted optimal value of $c$. Experiments show that such $c$ does produce a very good numerical solution to the PDE, even if the value of $c$ is not the experimentally optimal one. If one insists on finding the experimentally optimal value, it can be achieved by a moderate search. The stopping criterion is totally different from Kansa's approach and perhaps has never appeared in the literature. Our stopping criterion proves to be very reliable and does lead to the experimentally optimal value of $c$.

\section{Poisson equations}
The partial differential equations we deal with are of the form
\begin{align}
\left\{ \begin{array}{ll}
               u_{xx}(x,y)+u_{yy}(x,y)=f(x,y) & \mbox{for $(x,y)\in \Omega\backslash  \partial \Omega$}, \\
              u(x,y)=g(x,y)                              & \mbox{for $(x,y)\in \partial \Omega$}    
           \end{array}   
\right.
\end{align}
where $\Omega$ is the domain with boundary $\partial \Omega$, and $f,\ g$ are given functions. The reason we choose the Dirichlet condition as the boundary condition is that this setting is closer to function interpolation. For simplicity we let $\Omega$ be a square.

\subsection{1D experiment}
Although we are interested mainly in two-dimensional problems, as a prelude, a one-dimensional problem is illustrated and tested so that the reader can grasp the central idea and obtain a simple understanding for our approach.

In this experiment the solution function is $u(x)=e^{(-\sigma/2.1)x^{2}}$ where $\sigma=1$. It satisfies the equations
\begin{align}
\left\{ \begin{array}{ll}
               u_{xx}(x)=e^{(-1/2.1)x^{2}}[(-2/2.1)^{2}x^{2}-2/2.1],\\
              u(0)=1,\ u(10)=e^{-100/2.1}
             \end{array} 
   \right.  
\end{align}
in the domain $[0, 10]$. By Luh \cite{Lu3}, $u\in E_{\sigma},\ \sigma=1$, and the MN curves of Case 2. apply if we choose $\beta=-1$. The reason we adopt the inverse multiquadrics is that their programming is easier. We offer six MN curves in Figs. 1-6, which serve as the essential error bounds for the function interpolations. The number $b_{0}$, which greatly affects the MN curves, denotes the diameter of the interpolation domain. In these figures it is easily seen that as the fill distance $\delta$ decreases, i.e. the number of data points increases, the optimal values of $c$ move to 120 and are fixed there at last. Empirical results, as shown in \cite{Lu3}, show that one should choose $c=120$ to make the approximation.

\begin{figure}[h]
\centering
\includegraphics[scale=1.0]{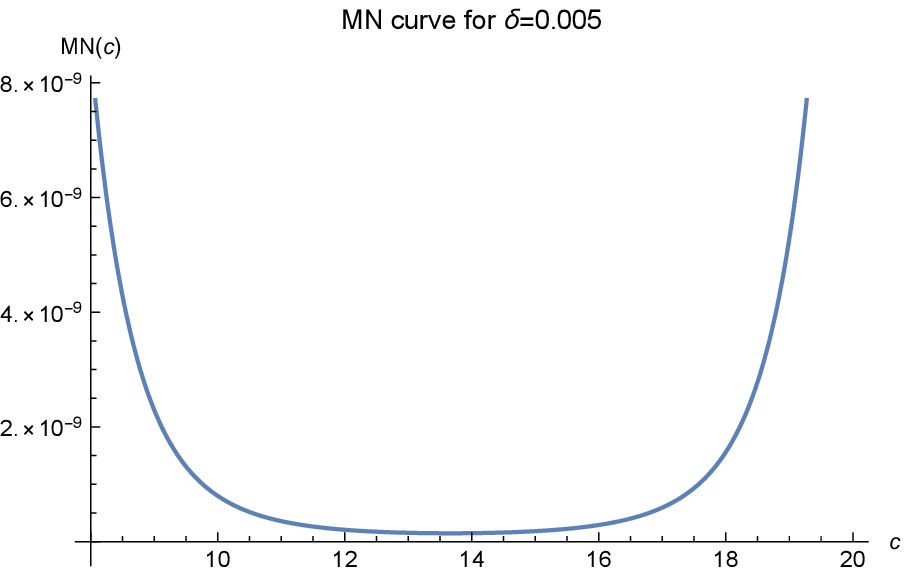}
\caption{Here $n=1,\ \beta=-1, b_{0}=10$ and $\sigma=1$.}

\includegraphics[scale=1.0]{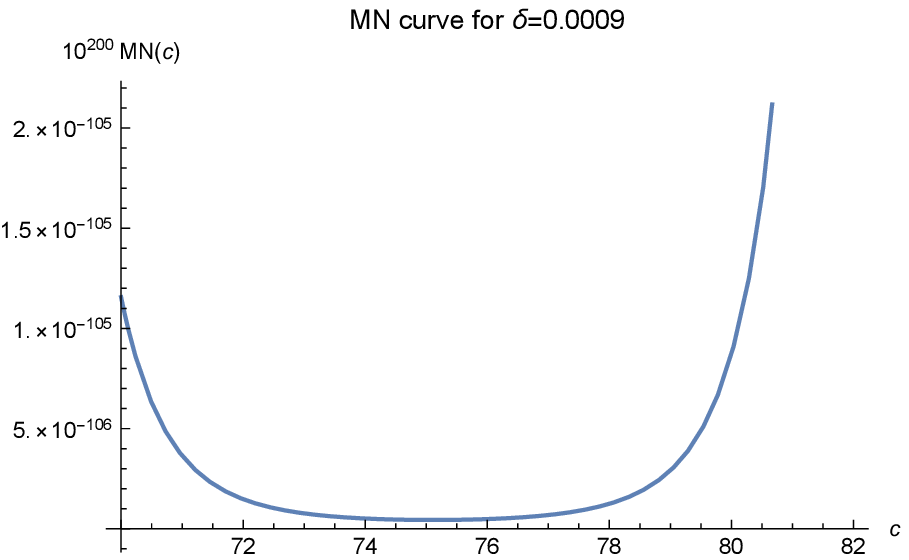}
\caption{Here $n=1,\ \beta=-1, b_{0}=10$ and $\sigma=1$.}

\includegraphics[scale=1.0]{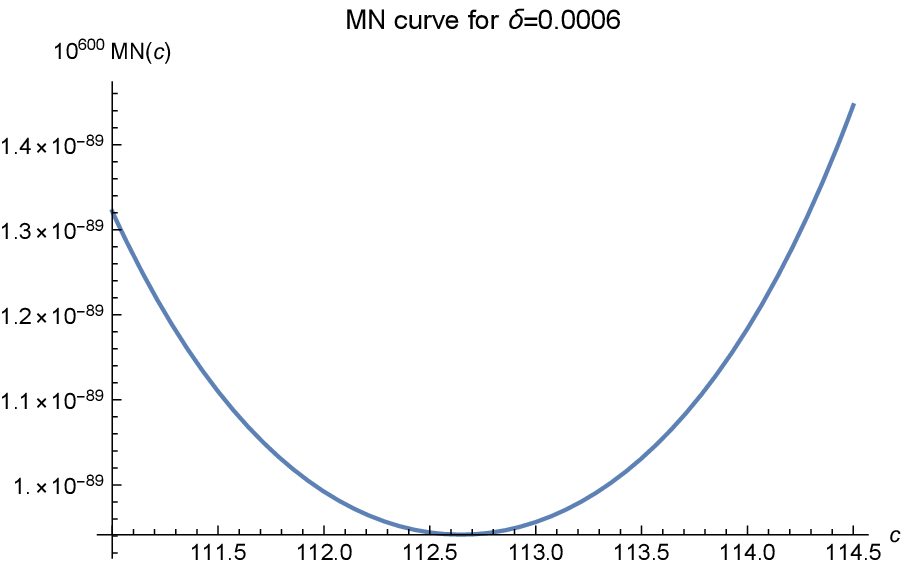}
\caption{Here $n=1,\ \beta=-1, b_{0}=10$ and $\sigma=1$.}

\end{figure}

\clearpage

\begin{figure}[h]
\centering
\includegraphics[scale=1.0]{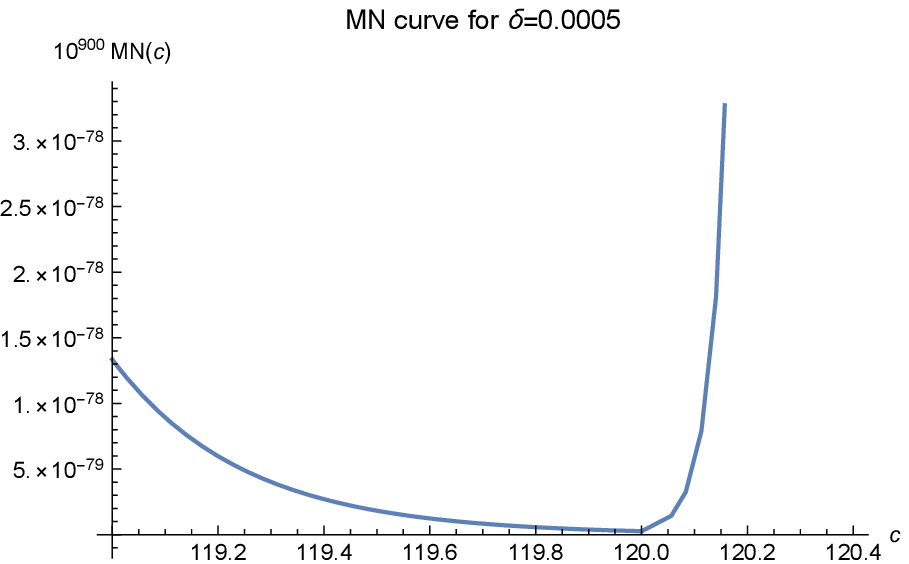}
\caption{Here $n=1,\ \beta=-1, b_{0}=10$ and $\sigma=1$.}

\includegraphics[scale=1.0]{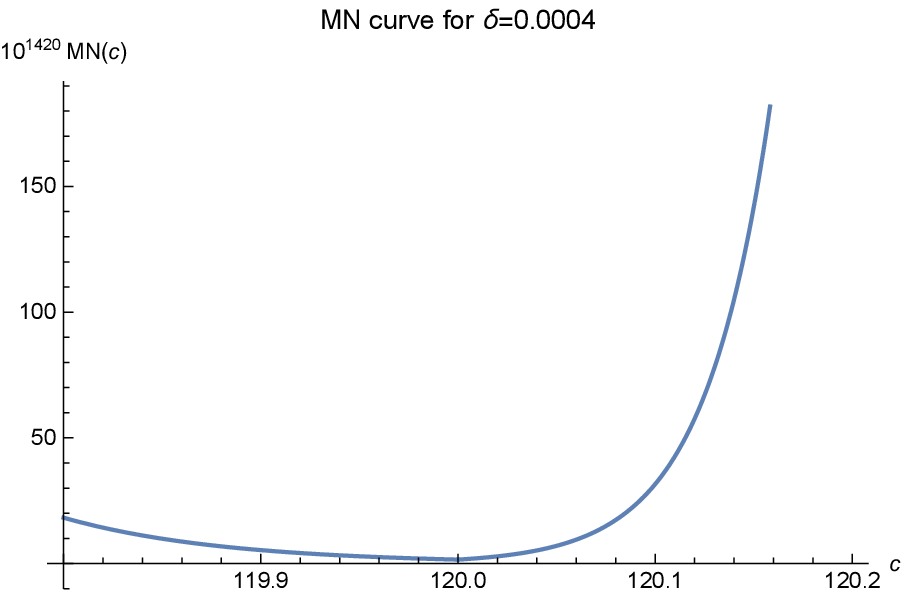}
\caption{Here $n=1,\ \beta=-1, b_{0}=10$ and $\sigma=1$.}

\includegraphics[scale=1.0]{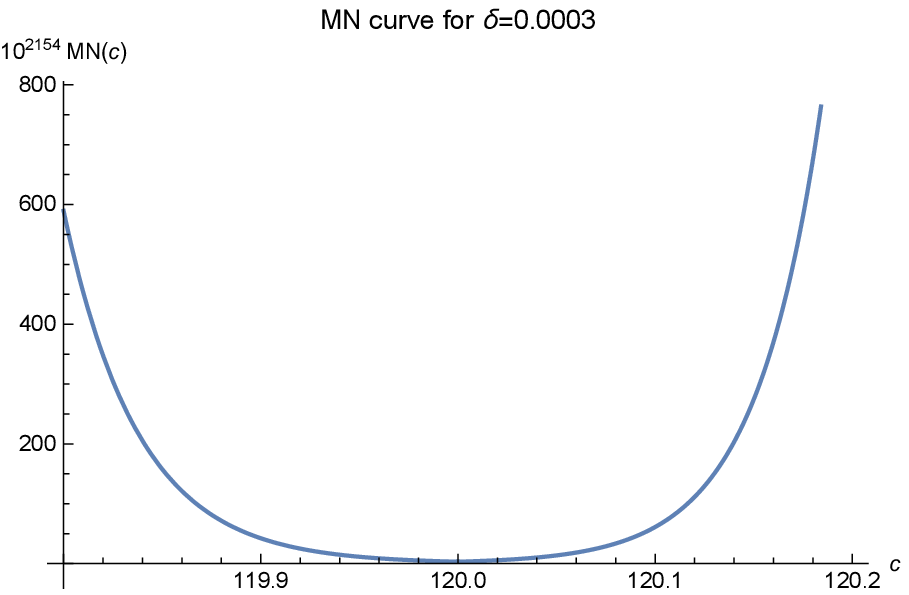}
\caption{Here $n=1,\ \beta=-1, b_{0}=10$ and $\sigma=1$.}

\end{figure}

\clearpage

Now we let $$\hat{u}(x)=\sum_{j=1}^{201}\lambda_{j}\phi(x-x_{j}),$$ where $\phi(x)=1/\sqrt{120^2+|x|^{2}}$, and require that $\hat{u}$ satisfy
$$ \hat{u}_{xx}(x_{j})=e^{(-1/2.1)x_{j}^{2}}[(-2/2.1)^{2}x_{j}^{2}-2/2.1]\ \mbox{for}\ j=2,\cdots,200,$$
and
$$\hat{u}(0)=1,\ \hat{u}(10)=e^{-100/2.1},$$
where $x_{1}=0,\ x_{201}=10,$ and $x_{j}=0.05(j-1)$ for $j=1,\cdots,201$.

This is a standard collocation setting. After solving the linear equations for $\lambda_{j}'s$, we tested $|u(x)-\hat{u}(x)|$ at 400 test points $z_{1},\cdots ,z_{400}$ evenly spaced in $[0,10]$ and found its root-mean-square error
$$RMS=\left\{\frac{1}{400}\sum_{i=1}^{400}|u(z_{i}-\hat{u}(z_{i})|^{2}\right\}^{1/2}=1.25\times 10^{-83}.$$ 
The condition number of the linear system is $4.4\times 10^{643}$. With the arbitrarily precise computer software Mathematica, we kept 800 effective digits to the right of the decimal point for each step of the calculation, successfully overcoming the problem of ill-conditioning. The computer time for solving the linear system was less than one second. We didn't test smaller fill distances $\delta's$ and different $c's$ because the RMS was already satisfactory.

\subsection{2D experiment}

Here the solution function is $u(x,y)=e^{-(\sigma/2.1)(x^{2}+y^{2})}$ where $\sigma=10^{-36}$. The domain is a large square with vertices $(0,0),\ (10^{16},0),\ (10^{16},10^{16})$ and $(0,10^{16})$. The function $u(x,y)$ satisfies
\begin{align}
u_{xx}(x,y)+u_{yy}(x,y)=-(2\sigma/2.1)e^{-(\sigma/2.1)(x^{2}+y^{2})}[2-(2\sigma/2.1)(x^{2}+y^{2})]
\end{align}
for $(x,y)$ in the interior of the domain $\Omega=\{(x,y)|\ 0\leq x\leq 10^{16},\ 0\leq y\leq 10^{16}\}$ and
\begin{align}
u(x,y)=e^{-(\sigma/2.1)(x^{2}+y^{2})}
\end{align}
for $(x,y)$ on the boundary $\partial \Omega$.

By Luh \cite{Lu3}, $u\in E_{\sigma}$ where $\sigma=10^{-36}$ and Case 1 of \cite{Lu3} applies. Five MN curves are shown in Figs. 7-11.

\begin{figure}[h]
\centering
\includegraphics[scale=1.0]{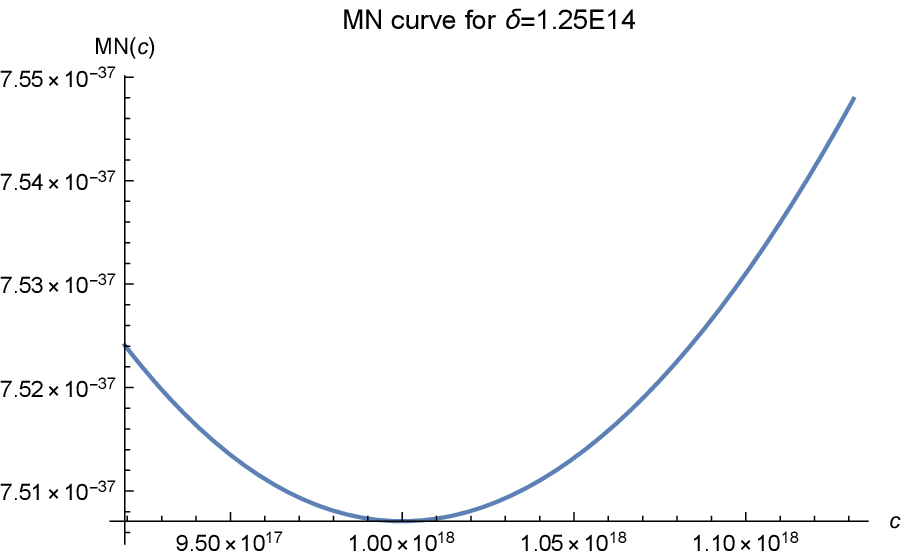}
\caption{Here $n=2,\ \beta=-1, b_{0}=\sqrt{2}$E$16$ and $\sigma=1$E$-36$.}

\includegraphics[scale=1.0]{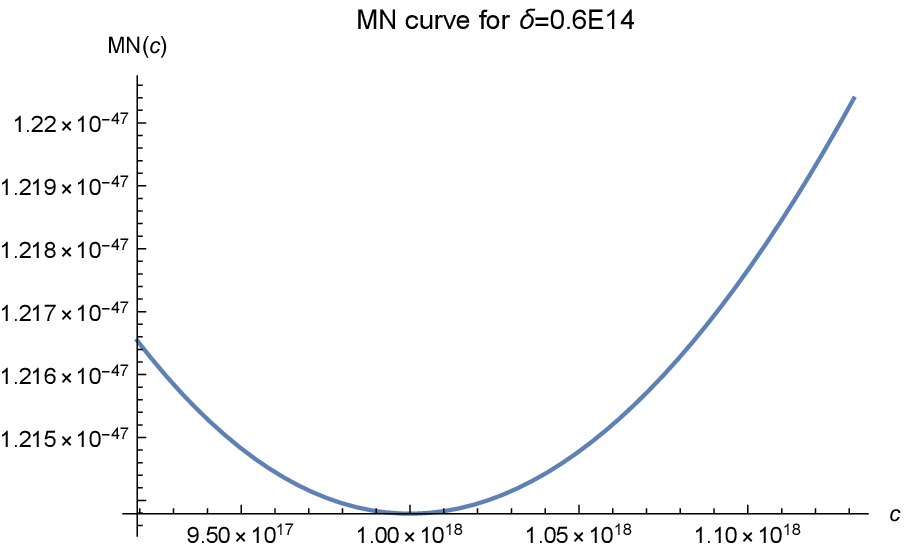}
\caption{Here $n=2,\ \beta=-1, b_{0}=\sqrt{2}$E$16$ and $\sigma=1$E$-36$.}

\includegraphics[scale=1.0]{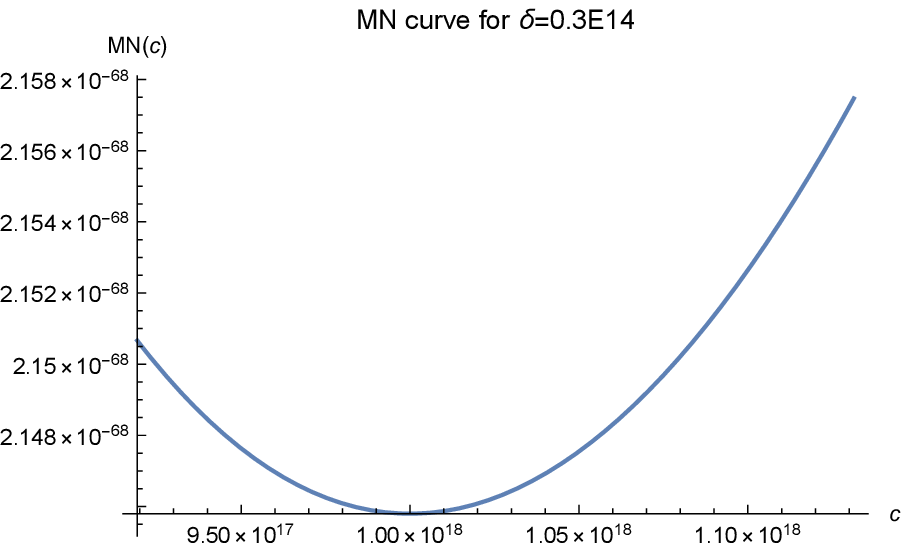}
\caption{Here $n=2,\ \beta=-1, b_{0}=\sqrt{2}$E$16$ and $\sigma=1$E$-36$.}

\end{figure}

\clearpage

\begin{figure}[h]
\centering
\includegraphics[scale=1.0]{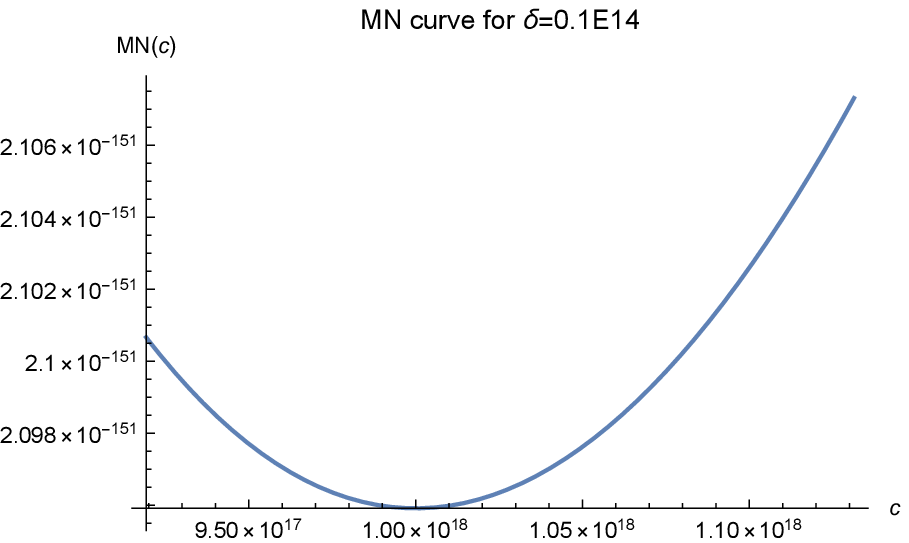}
\caption{Here $n=2,\ \beta=-1, b_{0}=\sqrt{2}$E$16$ and $\sigma=1$E$-36$.}

\includegraphics[scale=1.0]{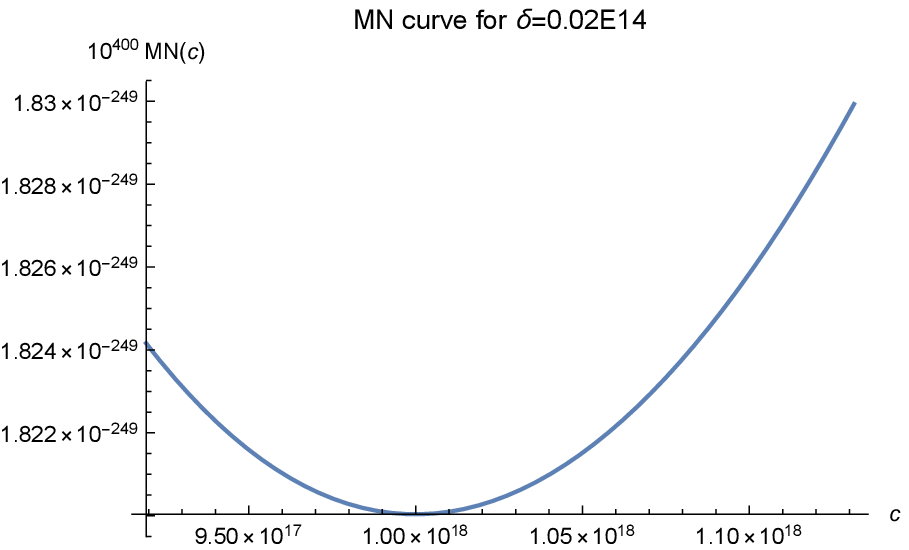}
\caption{Here $n=2,\ \beta=-1, b_{0}=\sqrt{2}$E$16$ and $\sigma=1$E$-36$.}

\end{figure}

All these curves show that one should choose $c=7000\cdot\sqrt{2}\cdot10^{14}\approx 0.99\cdot 10^{18}$ as the shape parameter in $\phi(x,y)=1/\sqrt{c^{2}+x^{2}+y^{2}}$. We let $\hat{u}(x,y):=\sum_{j=1}^{N_{d}}\lambda_{j}\phi(x-x_{j},y-y_{j})$ and require that it satisfy
$$\hat{u}_{xx}(x,y)+\hat{u}_{yy}(x,y)=-(2\sigma/2.1)e^{-(\sigma/2.1)(x^{2}+y^{2})}[2-(2\sigma/2.1)(x^{2}+y^{2})]$$
for the data points $(x,y)=(x_{j},y_{j})$ where $(x_{j},y_{j})$ belong to the interior of $\Omega$, i.e. $\Omega\backslash \partial \Omega$, and $N_{d}$ denotes the number of data points used. Also,
$$\hat{u}(x,y)=e^{-(\sigma/2.1)(x^{2}+y^{2})}$$
for $(x,y)=(x_{j},y_{j})$ where $(x_{j},y_{j})$ belong to the boundary $\partial \Omega$.

A grid of $41\times 41$ was adopted. Hence there are 1681 data points $(x_{j},y_{j})$ altogether. Among them 160 are boundary points where the Dirichlet condition occurs. Thus the fill distance is $\delta=1.25\sqrt{2}$E$14$. When applying the MN curves, we considered all the 1681 data points to be the interpolation points, even though it is not theoretically rigorous. As explained in Luh \cite{Lu3}, it is supposed to work well. However, something important must be pointed out. Although MN curves can be used to predict almost exactly the optimal value of $c$ for function interpolations, a moderate search may be needed if this approach is used in a non-rigorous way. We began with the theoretically predicted optimal value $c=7000\sqrt{2}\cdot10^{14}$, and tested two values nearby, one larger and the other smaller. Then we checked the RMS on the boundary for each $c$ and chose the direction which made the RMS smaller. Continuing choosing $c$ in this direction, we stopped when the RMS's began to grow. Our experiment shows that not many steps are needed, and the finally obtained $c$ does produce the best result.

The experimental results are presented in Table 1. Here $RMS,\ N_{d},\ N_{t},\ COND$ denote the root-mean-square error, number of data points, number of test points, and the condition number of the linear system, respectively. We use $RMSbdy$ to denote the root-mean-square-error of the approximation on the boundary, generated by 800 test points located on the boundary. In the entire domain $\Omega$, 6400 test points were used to generate the $RMS$'s. The most time-consuming command of solving the system of linear equations took about 30 minutes for each $c$. Although we adopted 1200 effective digits for each step of the calculation, it still worked with acceptable time efficiency.

\begin{table}[h]
\caption{$\delta=\sqrt{2}\times 1.25\times 10^{14},\ b_{0}=\sqrt{2}\times 10^{16},\ N_{d}=1681,\ N_{t}=6400$}
\centering
\tiny
\begin{tabular}{c lllll}\\[2ex]
\hline\hline \\ [1ex]
\large $c$ & \large $300 \sqrt{2}\cdot 10^{14}$ & \large $400 \sqrt{2}\cdot 10^{14}$ & \large$500 \sqrt{2}\cdot 10^{14}$ & \large $600 \sqrt{2}\cdot 10^{14}$ & \large$700 \sqrt{2}\cdot 10^{14}$  \\ [1ex]
\hline\\[1ex]

\large $RMS$ & \normalsize $ 5.2\cdot 10^{-134} $ & \normalsize $3.7\cdot 10^{-139} $ & \normalsize $ 1.5\cdot 10^{-143} $ & \normalsize $ 2.1\cdot 10^{-145} $ & \normalsize $5.8\cdot 10^{-147} $     \\ [1ex]
\hline\\[1ex]

\large $COND$ & \normalsize $1.1\cdot 10^{505}$ & \normalsize $5.8\cdot 10^{524}$ & \normalsize $1.2 \cdot 10^{540}$ & \normalsize $3.9 \cdot 10^{552}$ & \normalsize $1.5\cdot 10^{563}$  \\ [1ex]
\hline\\[1ex]

\large $RMSbdy$ & \normalsize $1.07\cdot 10^{-143}$ & \normalsize $7.2\cdot 10^{-149}$  & \normalsize $3.2\cdot 10^{-153}$ & \normalsize $4.1\cdot 10^{-155}$  & \normalsize $1.1\cdot 10^{-156}$\\[1ex]

\hline\hline \\[1ex]

\large $c$ & \large $800\sqrt{2}\cdot 10^{14}$ & \large $900\sqrt{2}\cdot 10^{14}$  & \normalsize $1000\sqrt{2}\cdot 10^{14}$ & \normalsize $1100\sqrt{2}\cdot 10^{14}$ & \normalsize $1200\sqrt{2}\cdot 10^{14}$ \\ [1ex]
\hline\\[1ex]

\large $RMS$ & \normalsize $ 3.6\cdot 10^{-148} $ & \normalsize $7.9\cdot 10^{-148} $ & \normalsize $8.2\cdot 10^{-146}$ & \normalsize $5.2\cdot 10^{-144}$ & \normalsize $6.0\cdot 10^{-143}$  \\ [1ex]
\hline \\[1ex]

\large $COND$ & \normalsize $2.1\cdot 10^{572}$ & \normalsize $2.6\cdot 10^{580}$ & \normalsize $4.4\cdot 10^{587}$ & \normalsize $1.5\cdot 10^{594}$ & \normalsize $1.4\cdot 10^{600}$  \\ [1ex]
\hline \\[1ex]

\large $RMSbdy$ & \normalsize $7.2\cdot 10^{-158}$  & \normalsize $1.7\cdot 10^{-157}$ & \normalsize $1.6\cdot 10^{-155}$ & \normalsize $1.0\cdot 10^{-153}$ & \normalsize $1.2\cdot 10^{-152}$   \\[1ex]

\hline\hline \\[1ex]

\large $c$ & \large $3000\sqrt{2}\cdot 10^{14}$ & \normalsize $5000\sqrt{2}\cdot 10^{14}$ & \normalsize $7000\sqrt{2}\cdot 10^{14}$ & \normalsize $8000\sqrt{2}\cdot 10^{14}$ \\[1ex]
\hline\\[1ex]

\large $RMS$ & \normalsize  $2.3\cdot 10^{-133}$ & \normalsize $7.1\cdot 10^{-127}$ & \normalsize $1.1\cdot 10^{-122}$ & \normalsize $5.4\cdot 10^{-121}$\\[1ex]
\hline \\[1ex]

\large $COND$ & \normalsize $1.1\cdot 10^{663}$ & \normalsize $1.2\cdot 10^{698}$ & \normalsize $1.5\cdot 10^{721}$  & \normalsize $2.2\cdot 10^{730}$\\[1ex]
\hline \\[1ex]

\large $RMSbdy$ & \normalsize  $4.4\cdot 10^{-143}$ & \normalsize $1.8\cdot 10^{-136}$ & \normalsize $5.1\cdot 10^{-132}$ & \normalsize $3.2\cdot 10^{-130}$\\[1ex]
\hline \\[1ex]

\end{tabular}

\label{001}
\end{table}
Note that the optimal value of $c$ is $800\sqrt{2}\cdot 10^{14}$ which coincides with the value chosen by our stopping criterion based on $RMSbdy$.
Obviously we could have got better RMS by increasing the number of data points, whereas we didn't do so because the approximation was already quite good.

\section{Final conclusion}
In physics many numerical solutions to PDEs are not bad, but truly good solutions are rarely seen. E. Kansa invented the collocation method and opened a new route to solving them. The combination of the c-theory and collocation does produce very good results as shown in our experiments. Maybe this is just a starting point. We are still facing a huge challenge and have a lot of work to do in the future.

\end{document}